\documentclass[12pt,draft]{article}

\usepackage[cp1251]{inputenc}
\usepackage[T2A]{fontenc}
\usepackage[english,ukrainian]{babel}
\usepackage{amsmath,amsfonts,amssymb}
\usepackage{geometry}

\begin{document}

\begin{center}

\large

\textbf{ELLIPTIC PROBLEMS \\WITH NONCLASSICAL BOUNDARY CONDITIONS \\IN AN EXTENDED SOBOLEV SCALE}

\medskip

\textbf{A.A. Murach, I.S. Chepurukhina}

\medskip

\normalsize

Institute of Mathematics of the National Academy of Sciences of Ukraine, Kyiv

\medskip

E-mail: murach@imath.kiev.ua, Chepurukhina@gmail.com

\bigskip\bigskip

\large

\textbf{ЕЛІПТИЧНІ ЗАДАЧІ \\З НЕКЛАСИЧНИМИ КРАЙОВИМИ УМОВАМИ\\ У РОЗШИРЕНІЙ СОБОЛЄВСЬКІЙ ШКАЛІ}\footnote{Публікація містить результати досліджень, проведених за грантом Президента України за конкурсним проектом Ф82/45932 Національного фонду досліджень України.}

\medskip

\textbf{О.О. Мурач, І.С. Чепурухіна}

\normalsize

\end{center}

\medskip

\noindent We consider elliptic problems with nonclassical boundary conditions that contain additional unknown functions on the border of the domain of the elliptic equation and also contain boundary operators of  higher orders with respect to the order of this equation. We investigate the solvability of the indicated problems and properties of their solutions in an extended Sobolev scale. It consists of Hilbert generalized Sobolev spaces for which the order of regularity is a general radial function $\mathrm{RO}$-varying in the sense of Avakumovi\'{c} at infinity. We establish a theorem on the Fredholm property of the indicated problems on appropriate pairs of these spaces and theorems on the regularity and \textit{a~priori} estimate of the generalized solutions to the problems. We obtain exact sufficient conditions for components of these solutions to be continuously differentiable.

\medskip



\bigskip\bigskip

\noindent \textbf{Вступ.} Сучасна теорія еліптичних крайових задач найбільш детально розроблена у випадку регулярних крайових умов (див., наприклад, [\ref{LionsMagenes71}, розд. 2]). У цьому випадку виконується класична формула Гріна, що суттєво полегшує дослідження задач, особливо тоді, коли праві частини є узагальненими функціями. Втім, у застосуваннях виникають еліптичні задачі з більш загальними, некласичними, крайовими умовами, для яких не виконується згадана формула Гріна. Зокрема, так буває, коли порядки крайових операторів більші за порядок еліптичного рівняння, або рівні йому. До таких задач належать і крайові задачі з додатковими невідомими функціями на межі евклідової області, де задано еліптичне рівняння; їх уперше розглянув Б. Лаврук [\ref{Lawruk63a}].

У нашій роботі досліджуються еліптичні задачі з некласичними крайовими умовами, які містять додаткові невідомі функції на межі області та   крайові оператори порядків вищих, ніж порядок еліптичного рівняння. Такі задачі вивчено у просторах Соболєва в [\ref{KozlovMazyaRossmann97}, п. 4.1] і [\ref{Roitberg99}, п.~2.4]. Мета нашої роботи --- встановити теореми про характер розв'язності і властивості розв'язків цих задач в узагальнених просторах Соболєва, які утворюють розширену соболєвську шкалу [\ref{MikhailetsMurach13UMJ3}; \ref{MikhailetsMurach14}, п. 2.4.2]. Показником регулярності для цих просторів служить не число (як для просторів Соболєва), а радіальна функція, RO-змінна на нескінченності. Використання функціонального параметра замість числового дозволяє більш тонко охарактеризувати регулярність узагальнених розв'язків досліджуваних задач, зокрема, отримати точні достатні умови неперервної диференційовності компонент цих розв'язків.

\textbf{1. Постановка задачі.} Нехай $\Omega$~--- довільна обмежена область у евклідовому просторі $\mathbf{R}^{n}$, де $n\geq2$. Припустимо, що її межа $\Gamma$ є нескінченно гладким замкненим многовидом вимірності $n-1$, причому $C^{\infty}$- структура на $\Gamma$ породжена $\mathbf{R}^{n}$.

Нехай задано цілі числа $q\geq1$, $\lambda\geq1$, $m_{1},\ldots,m_{q+\lambda}$ та $r_{1},\ldots,r_{\lambda}$.
В~області $\Omega$ розглядаємо крайову задачу вигляду
\begin{gather}\label{13f1}
Au=f\quad\mbox{в}\quad\Omega,\\
B_{j}u+\sum_{k=1}^{\lambda}C_{j,k}v_{k}=g_{j}\quad\mbox{на}\quad\Gamma,
\quad j=1,...,q+\lambda.\label{13f2}
\end{gather}
Тут $A:=A(x,D)$~--- лінійний диференціальний оператор (л.д.о.) на $\overline{\Omega}:=\Omega\cup\Gamma$ парного порядку $2q$, кожне $B_{j}:=B_{j}(x,D)$~--- крайовий л.д.о. на $\Gamma$ порядку $\mathrm{ord}\,B_{j}\leq m_{j}$, а кожне $C_{j,k}:=C_{j,k}(x,D_{\tau})$~--- дотичний л.д.о. на $\Gamma$ порядку $\mathrm{ord}\,C_{j,k}\leq m_{j}+r_{k}$. (Як звичайно, диференціальні оператори від'ємного порядку вважаються нуль-операторами.) Усі коефіцієнти цих л.д.о. є нескінченно гладкими функціями, заданими на $\overline{\Omega}$ і $\Gamma$ відповідно. Функція $u$ на $\Omega$ і функції  $v_{1},\ldots,v_{\lambda}$ на $\Gamma$ є шуканими у крайовій задачі \eqref{13f1}, \eqref{13f2}. У~статті всі функції та розподіли вважаються комплекснозначними і тому розглядаються комплексні лінійні функціональні простори.

Покладемо $m:=\max\{m_1,\ldots,m_{q+\lambda}\}$. У роботі розглядаємо некласичний випадок, коли $m\geq2q$. Вважаємо також, що $m\geq-r_{k}$ для кожного цілого $k\in[1,\lambda]$ (якщо $m+r_{k}<0$ для деякого вказаного $k$, то усі оператори $C_{1,k}$,..., $C_{q+\lambda,k}$ дорівнюють нулю і тому шукана функція $v_{k}$ відсутня у крайових умовах \eqref{13f2}).

Надалі припускаємо, що крайова задача \eqref{13f1}, \eqref{13f2} є еліптичною в області $\Omega$, тобто л.д.о. $A$ правильно еліптичний на $\overline{\Omega}$, а система крайових умов \eqref{13f2} накриває $A$ на $\Gamma$ (див., наприклад, [\ref{KozlovMazyaRossmann97}, п.~3.1.2]).

Прикладом еліптичної крайової задачі вигляду \eqref{13f1}, \eqref{13f2} є така задача:
\begin{gather*}
\Delta u=f\quad\mbox{в}\quad\Omega,\\
\partial_{\nu}^{p}u+iav_1=g_{1},\quad
\partial_{\nu}^{p+2l}u+b\Delta_{\Gamma}^{l}v_1=g_{2}
\quad\mbox{на}\quad\Gamma.
\end{gather*}
Тут довільно вибрано цілі числа $p\geq0$ і $l\geq1$ та дійсні функції  $a,b\in C^{\infty}(\Gamma)$ такі, що $|a(x)|+|b(x)|\neq0$ для довільного $x\in\Gamma$. Як звичайно, $\Delta$~--- оператор Лапласа в $\mathbf{R}^{n}$ і $\Delta_{\Gamma}$~--- оператор Бельтрамі\,--\,Лапласа на $\Gamma$, а $\partial_{\nu}$~--- оператор диференціювання вздовж внутрішньої нормалі до межі $\Gamma$ області $\Omega$. Для цієї задачі $m=p+2l\geq\mathrm{ord}\,\Delta$. Якщо $a(x_0)=0$ для деякого  $x_0\in\Gamma$, то не можна вилучити невідому функцію $v_1$ з крайових умов і зберегти при цьому гладкість коефіцієнтів і правих частин крайових умов, що отримуються.

\textbf{2. Розширена соболєвська шкала} утворена гільбертовими узагальненими соболєвськими просторами $H^{\alpha}$, для яких показником регулярності служить довільна функція $\alpha$ з класу RO. За означенням, він складається з усіх вимірних за Борелем функцій $\alpha:\nobreak[1,\infty)\rightarrow(0,\infty)$, для яких існують числа $b>1$ і $c\geq1$ такі, що $c^{-1}\leq\alpha(\lambda t)/\alpha(t)\leq c$ для довільних $t\geq1$ і $\lambda\in[1,b]$ (сталі $b$ і $c$ можуть залежати від $\alpha$). Ці функції називають RO-змінними за В.Г. Авакумовичем на нескінченності (див., наприклад, [\ref{Seneta85}, с. 86]).

Клас RO має простий опис [\ref{Seneta85}, с.~87]:
$$
\alpha\in\mathrm{RO}\;\;\Longleftrightarrow\;\;
\alpha(t)=\exp\Biggl(\beta(t)+
\int\limits_{1}^{\:t}\frac{\gamma(\tau)}{\tau}\;d\tau\Biggr)\;\,
\mbox{для}\;\,t\geq1,
$$
де дійсні функції $\beta$ і $\gamma$ вимірні за Борелем і обмежені на півосі $[1,\infty)$. Як відомо [\ref{Seneta85}, с.~88], для кожної функції $\alpha\in\mathrm{RO}$ існують дійсні числа $s_{0}<s_{1}$ і $c_{0},c_{1}\in(0,\infty)$ такі, що
\begin{equation}\label{13f8}
c_{0}\lambda^{s_{0}}\leq\alpha(\lambda t)/\alpha (t)\leq
c_{1}\lambda^{s_{1}} \quad\mbox{для всіх}\quad t,\lambda\in[1,\infty).
\end{equation}
Позначимо через $\sigma_{0}(\alpha)$ супремум усіх чисел $s_{0}$, для яких виконується ліва нерівність в \eqref{13f8}, а через $\sigma_{1}(\alpha)$~--- інфімум усіх чисел $s_{1}$, для яких виконується права нерівність в \eqref{13f8}. Числа $\sigma_{0}(\alpha)$ і $\sigma_{1}(\alpha)$ називають відповідно нижнім і верхнім індексами Матушевської функції $\alpha$ (див. [\ref{BinghamGoldieTeugels89}, п.~2.1.2]).

Нехай $\alpha\in\mathrm{RO}$. За означенням, лінійний простір
$H^{\alpha}(\mathbf{R}^{n})$, де $n$~--- натуральне число, складається з усіх повільно зростаючих розподілів $w$ на $\mathbf{R}^{n}$ таких, що їх перетворення Фур'є $\widehat{w}$ локально інтегровне за Лебегом на $\mathbf{R}^{n}$ і задовольняє умову
\begin{equation*}
\|w\|_{\alpha,\mathbf{R}^{n}}^{2}:=\int\limits_{\mathbf{R}^{n}}
\alpha^2(\langle\xi\rangle)\,|\widehat{w}(\xi)|^2\,d\xi<\infty.
\end{equation*}
Тут $\langle\xi\rangle:=(1+|\xi|^{2})^{1/2}$ є згладжений модуль вектора $\xi\in\mathbf{R}^{n}$. За означенням, $\|w\|_{\alpha,\mathbf{R}^{n}}$~--- норма у цьому просторі.

Простір $H^{\alpha}(\mathbf{R}^{n})$ є гільбертів ізотропний випадок просторів, уведених і досліджених Л.~Хермандером [\ref{Hermander65}, п.~2.2] та Л.Р.~Волевичем і Б.П.~Панеяхом [\ref{VolevichPaneah65}, \S~2].
Якщо функція $\alpha(t)\equiv t^{s}$ степенева, то $H^{\alpha}(\mathbf{R}^{n})$ стає (гільбертовим) простором
Соболєва $H^{(s)}(\mathbf{R}^{n})$ дійсного порядку $s$. Узагалі,
\begin{equation}\label{13f9}
(s_{0}<\sigma_{0}(\alpha),\;\sigma_{1}(\alpha)<s_{1})\;\;\Longrightarrow
\;\;H^{(s_1)}(\mathbf{R}^{n})\subset H^{\alpha}(\mathbf{R}^{n})\subset
H^{(s_0)}(\mathbf{R}^{n}),
\end{equation}
де вкладення неперервні й щільні.

Клас функціональних просторів $H^{\alpha}(\mathbf{R}^{n})$, де $\alpha\in\mathrm{RO}$, виділений в [\ref{MikhailetsMurach13UMJ3}, п. 2] і названий розширеною соболєвською шкалою (р.с.ш.) на $\mathbf{R}^{n}$.
Нам потрібні її аналоги для області $\Omega\subset\mathbf{R}^{n}$ та її межі~$\Gamma$. Наведемо відповідні означення (див. [\ref{MikhailetsMurach15ResMath1}, п. 2] і [\ref{MikhailetsMurach14}, п. 2.4.2]). Тепер ціле $n\geq2$.

Нехай, як і раніше, $\alpha\in\mathrm{RO}$. За означенням, лінійний простір $H^{\alpha}(\Omega)$ складається зі звужень в область $\Omega$
усіх розподілів $w\in H^{\alpha}(\mathbf{R}^{n})$. Цей простір наділений нормою
$$
\|u\|_{\alpha,\Omega}:=\inf\bigl\{\,\|w\|_{\alpha,\mathbf{R}^{n}}:
w\in H^{\alpha}(\mathbf{R}^{n}),\;w=u\;\,\mbox{в}\;\,\Omega\,\bigr\},
$$
де $u\in H^{\alpha}(\Omega)$. Він є сепарабельним гільбертовим простором відносно цієї норми, причому множина $C^{\infty}(\overline{\Omega})$ щільна у ньому.

Лінійний простір $H^{\alpha}(\Gamma)$ складається, коротко кажучи, з усіх розподілів на $\Gamma$, які в локальних координатах дають елементи простору $H^{\alpha}(\mathbf{R}^{n-1})$. Дамо детальне означення. Виберемо з $C^{\infty}$-структури на $\Gamma$ деякий скінченний набір локальних карт $\pi_j: \mathbf{R}^{n-1}\leftrightarrow\Gamma_{j}$, де $j\in\{1,\ldots,p\}$, а відкриті множини $\Gamma_{1}$,..., $\Gamma_{p}$ утворюють покриття многовиду $\Gamma$. Окрім того, нехай функції $\chi_j\in C^{\infty}(\Gamma)$, де $j\in\{1,\ldots,p\}$, утворюють розбиття одиниці на $\Gamma$, яке задовольняє умову $\mathrm{supp}\,\chi_j\subset \Gamma_j$. За означенням, лінійний простір $H^{\alpha}(\Gamma)$ складається з усіх розподілів $h$ на $\Gamma$ таких, що $(\chi_{j}h)\circ\pi_{j}\in H^{\alpha}(\mathbf{R}^{n-1})$ для кожного $j\in\{1,\ldots,p\}$; тут $(\chi_{j}h)\circ\pi_{j}$ є зображенням розподілу $\chi_{j}h$ у локальній карті $\pi_{j}$. Простір $H^{\alpha}(\Gamma)$ наділений нормою
$$
\|h\|_{\alpha,\Gamma}:=
\biggl(\,\sum_{j=1}^{p}\,\|(\chi_{j}h)\circ\pi_{j}\|
_{\alpha,\mathbf{R}^{n-1}}^{2}\biggr)^{1/2}.
$$
Він гільбертів і сепарабельний відносно цієї норми та з точністю до еквівалентності норм не залежить від вибору атласу і розбиття одиниці [\ref{MikhailetsMurach14}, теорема~2.21]. Множина $C^{\infty}(\Gamma)$ щільна у цьому просторі.

Гільбертові простори $H^{\alpha}(G)$, де $\alpha\in\mathrm{RO}$, утворюють р.с.ш. на $G\in\{\Omega,\Gamma\}$. Якщо $\alpha(t)\equiv t^{s}$ для деякого дійсного числа $s$, то $H^{\alpha}(G)$ є соболєвським простором $H^{(s)}(G)$ порядку~$s$. Для р.с.ш. на $G$ зберігається властивість \eqref{13f9}, якщо у ній замінити $\mathbf{R}^{n}$ на $G$; при цьому обидва вкладення просторів будуть компактними і щільними.

Розширена соболєвська шкала має важливі інтерполяційні властивості, які відіграють ключову роль у її застосуваннях, зокрема до еліптичних операторів і еліптичних крайових задач. Вона отримується інтерполяцією з функціональним параметром пар гільбертових просторів Соболєва, складається з усіх гільбертових просторів, інтерполяційних відносно цих пар і замкнена відносно інтерполяції з функціональним параметром пар гільбертових просторів (див. [\ref{MikhailetsMurach14}, п. 2.4.2] і
[\ref{MikhailetsMurach15ResMath1}, п. 2 і 5]).

\textbf{3. Основні результати} статті стосуються властивостей еліптичної крайової задачі \eqref{13f1}, \eqref{13f2} у просторах, що належать до р.с.ш.

\textbf{Теорема 1.} \it Розглянемо відображення $(u,v_{1},...,v_{\lambda})\to(f,g_{1},...,g_{q+\lambda})$, де $u\in C^{\infty}(\overline{\Omega})$ і $v_{1},\ldots,v_{\lambda}\in C^{\infty}(\Gamma)$, а функції $f$ і $g_{1}$,...,$g_{q+\lambda}$ означені за формулами \eqref{13f1} і \eqref{13f2}. Воно продовжується єдиним чином (за неперервністю) до обмеженого лінійного оператора
\begin{equation}\label{13f12}
\Lambda:H^{\eta}(\Omega)\oplus
\bigoplus_{k=1}^{\lambda}H^{\eta\rho^{r_{k}-1/2}}(\Gamma)\to
H^{\eta\rho^{-2q}}(\Omega)\oplus
\bigoplus_{j=1}^{q+\lambda}H^{\eta\rho^{-m_{j}-1/2}}(\Gamma)
\end{equation}
для довільного параметра $\eta\in\mathrm{RO}$ такого, що $\sigma_{0}(\eta)>m+1/2$. Цей оператор нетерів. Його скінченновимірне ядро лежить у просторі $C^{\infty}(\overline{\Omega})\times(C^{\infty}(\Gamma))^{\lambda}$ і разом зі скінченним індексом не залежить від $\eta$. \rm

Тут і далі використано функціональний параметр $\rho(t):=t$ аргументу $t\geq1$ для того, щоб не записувати цей аргумент
у позначеннях узагальнених соболєвських просторів. Отже, параметр $\eta\rho^{l}$, де $l$~--- дійсне число, позначає функцію $\eta(t)t^{l}$. Звісно, $\eta\in\mathrm{RO}\Leftrightarrow\eta\rho^{l}\in\mathrm{RO}$, причому
$\sigma_{j}(\eta\rho^{l})=\sigma_{j}(\eta)+l$, де $j\in\{0,1\}$.

Умова $\sigma_0(\eta)>m+1/2$ суттєва у теоремі~1. Справді, якщо $\mathrm{ord}\,B_{j}=m$ для деякого цілого $j\in[1,q+\lambda]$, то відображення $u\mapsto B_{j}u$, де $u\in C^{\infty}(\overline{\Omega})$,
не можна продовжити до неперервного лінійного оператора, що діє з соболєвського простору $H^{(m+1/2)}(\Omega)$ у лінійний топологічний простір усіх розподілів на~$\Gamma$.

Дослідимо властивості узагальнених розв'язків крайової задачі \eqref{13f1}, \eqref{13f2}. Позначимо через $D^{\eta}$ область визначення оператора \eqref{13f12}, а через $D^{m+1/2+}$ об'єднання усіх просторів $D^{\eta}$, де параметр $\eta\in\mathrm{RO}$ задовольняє умову $\sigma_{0}(\eta)>m+1/2$. Вектор
\begin{equation}\label{generalized-solution}
(u,v):=(u,v_{1},\ldots,v_{\lambda})\in\mathcal{D}^{m+1/2+}
\end{equation}
називаємо узагальненим розв'язком крайової задачі \eqref{13f1}, \eqref{13f2} з правою частиною
$$
(f,g):=(f,g_{1},\ldots,g_{q+\lambda})\in D'(\Omega)\times\bigl(D'(\Gamma)\bigr)^{q+\lambda},
$$
якщо $\Lambda(u,v)=(f,g)$, де $\Lambda$ є оператор \eqref{13f12} для деякого параметра $\eta\in \mathrm{RO}$ такого, що $\sigma_0(\eta)>m+1/2$. Це означення коректне, оскільки воно не залежить від $\eta$. Тут, як звичайно, $D'(\Omega)$ і $D'(\Gamma)$ позначають лінійні топологічні простори усіх розподілів в $\Omega$ і на $\Gamma$ відповідно.

Нехай відкрита множина $V\subset\mathbf{R}^{n}$ така, що $\Omega_{0}:=\Omega\cap V\neq\varnothing$. Покладемо $\Gamma_{0}:=\Gamma\cap V$ (можливий випадок, коли $\Gamma_{0}=\varnothing$). Позначимо через $H^{\alpha}_{\mathrm{loc}}(\Omega_{0},\Gamma_{0})$, де $\alpha\in\mathrm{RO}$, лінійний простір усіх розподілів $u\in D'(\Omega)$ таких, що $\chi u\in H^{\alpha}(\Omega)$ для кожної функції $\chi\in C^{\infty}(\overline{\Omega})$, носій якої задовольняє умову $\mathrm{supp}\,\chi\subset\Omega_0\cup\Gamma_{0}$. Аналогічно, позначимо через $H^{\alpha}_{\mathrm{loc}}(\Gamma_{0})$ лінійний простір усіх розподілів $h\in D'(\Gamma)$ таких, що $\chi h\in H^{\alpha}(\Gamma)$ для кожної функції $\chi\in C^{\infty}(\Gamma)$, яка задовольняє умову $\mathrm{supp}\,\chi\subset\Gamma_{0}$.

\textbf{Теорема 2.} \it Нехай $\eta\in\mathrm{RO}$ і $\sigma_0(\eta)>m+1/2$. Припустимо, що  вектор \eqref{generalized-solution} є узагальненим розв'язком еліптичної крайової задачі \eqref{13f1}, \eqref{13f2}, праві частини якої задовольняють умови $f\in H^{\eta\rho^{-2q}}_{\mathrm{loc}}(\Omega_{0},\Gamma_{0})$ і $g_{j}\in H^{\eta\rho^{-m_{j}-1/2}}_{\mathrm{loc}}(\Gamma_{0})$ для кожного цілого $j\in[1,q+\lambda]$. Тоді цей розв'язок має такі властивості: $u\in H^{\eta}_{\mathrm{loc}}(\Omega_{0},\Gamma_{0})$ і $v_{k}\in H^{\eta\rho^{r_{k}-1/2}}_{\mathrm{loc}}(\Gamma_{0})$ для кожного цілого $k\in[1,\lambda]$. \rm

Якщо $\Omega_{0}=\Omega$ і $\Gamma_{0}=\Gamma$, то $H^{\alpha}_{\mathrm{loc}}(\Omega_{0},\Gamma_{0})=H^{\alpha}(\Omega)$ і $H^{\alpha}_{\mathrm{loc}}(\Gamma_{0})=H^{\alpha}(\Gamma)$. У цьому випадку теорема 2 стосується глобальної регулярності правих частин і розв'язку досліджуваної крайової задачі.

Цю теорему доповнює апріорна оцінка узагальненого розв'язку \eqref{generalized-solution}. Позначимо через $\|\cdot\|'_{\eta}$ і $\|\cdot\|''_{\eta\rho^{-2q}}$ норми у гільбертових просторах, на парі яких діє оператор \eqref{13f12}.

\textbf{Теорема 3.} \it Нехай $\eta\in\mathrm{RO}$ і $\sigma_0(\eta)>m+1/2$. Припустимо, що  вектор \eqref{generalized-solution} задовольняє умову теореми $2$. Довільно виберемо число $l>0$ і функції $\chi,\zeta\in C^{\infty}(\overline{\Omega})$ такі, що $\mathrm{supp}\,\chi\subset\mathrm{supp}\,\zeta\subset
\Omega_{0}\cup\Gamma_{0}$ та $\zeta=1$ у деякому околі множини $\mathrm{supp}\,\chi$. Тоді
\begin{equation*}
\|\chi(u,v)\|'_{\eta}\leq c\,\bigl(\,\|\zeta(f,g)\|''_{\eta\rho^{-2q}}+
\|\zeta(u,v)\|'_{\eta\rho^{-l}}\bigl)
\end{equation*}
для деякого числа $c>0$, незалежного від $(u,v)$ і $(f,g)$. \rm

У припущенні про те, що функція $\eta$ є правильно змінною за Й.Караматою на нескінченності, версії теорем 1--3 встановлено недавно в [\ref{ChepuruhinaKasirenko17}, п. 4]. Нагадаємо [\ref{Seneta85}, п. 1.1], що воно означає таке: $\eta(pt)/\eta(t)\to p^{\sigma}$ при $t\to\infty$ для кожного додатного числа $p>0$ і деякого дійсного числа $\sigma$; у цьому випадку індекси Матушевської функції  $\eta$ дорівнюють $\sigma$.
Випадки, коли $\lambda=0$ або усі $m_{j}\leq2q-1$ досліджені відповідно в [\ref{KasirenkoMurachUMJ3}] і [\ref{Chepuruhina15Coll2}]. Для соболєвських просторів ці теореми або їх версії відомі, див., наприклад,  монографії [\ref{KozlovMazyaRossmann97}, п. 4.1] і [\ref{Roitberg99}, п.~2.4].

\textbf{4. Застосування.} За допомогою р.с.ш. і теореми вкладення Л. Хермандера [\ref{Hermander65}, теорема 2.2.7] отримуємо тонкі й точні достатні умови $l$ разів неперервної диференційовності компонент узагальненого розв'язку \eqref{generalized-solution} еліптичної крайової задачі \eqref{13f1}, \eqref{13f2}. Нижче множини $\Omega_{0}$ і $\Gamma_{0}$ є такими як у п. 3, а $C^{l}(G)$ --- простір $l$ разів неперервно диференційовних функцій на множині $G\subset\mathbf{R}^{n}$.

\textbf{Теорема 4.} \it Нехай додатне ціле число $l>m+1/2-n/2$. Припустимо, що вектор \eqref{generalized-solution} задовольняє умову теореми $2$ для деякого функціонального параметра $\eta\in\mathrm{RO}$ такого, що $\sigma_{0}(\eta)>m+1/2$ і
\begin{equation}\label{13f19}
\int\limits_1^{\infty}t^{2l+n-1}\eta^{-2}(t)dt<\infty.
\end{equation}
Тоді $u\in C^{l}(\Omega_{0}\cup\Gamma_{0})$. Умова \eqref{13f19} є точною.
\rm

У цій теоремі умова $l>m+1/2-n/2$ є природною з огляду на включення \eqref{generalized-solution}, оскільки воно тягне за собою властивість $u\in C^{l}(\overline{\Omega})$ для будь-якого додатного цілого числа  $l\leq m+1/2-n/2$.

\textbf{Теорема 5.} \it Нехай $k\in\{1,...,\lambda\}$, а додатне ціле число $l>m+r_{k}+1/2-n/2$. Припустимо, що вектор \eqref{generalized-solution} задовольняє умову теореми $2$ для деякого функціонального параметра $\eta\in\mathrm{RO}$ такого, що $\sigma_{0}(\eta)>m+1/2$ і
\begin{equation}\label{13f20}
\int\limits_1^{\infty}t^{2(l-r_{k})+n-1}\eta^{-2}(t)dt<\infty.
\end{equation}
Тоді $v_{k}\in C^{l}(\Gamma_{0})$. Умова \eqref{13f20} є точною.
\rm

У цій теоремі умова на $l$ природна з огляду на включення \eqref{generalized-solution}, бо воно тягне за собою властивість $v_{k}\in C^{l}(\overline{\Omega})$ для кожного додатного цілого числа  $l\leq m+r_{k}+1/2-n/2$.

За допомогою цих теорем отримуємо достатню умову класичності узагальненого розв'язку \eqref{generalized-solution} досліджуваної задачі, тобто умову, за якої $u\in C^{2q}(\Omega)\cap C^{m}(U_{\delta}\cup\Gamma)$ для деякого числа $\delta>0$, а $v_{k}\in C^{m+r_{k}}(\Gamma)$ для кожного $k\in\{1,\ldots,\lambda\}$. Тут покладаємо $U_{\delta}:=\{x\in\Omega:\mathrm{dist}(x,\Gamma)<\delta\}$. Якщо розв'язок задачі класичний, то її ліві частини обчислюються за допомогою класичних похідних і є неперервними функціями на $\Omega$ і $\Gamma$ відповідно.

\textbf{Теорема 6.} \it Припустимо, що вектор \eqref{generalized-solution} є узагальненим розв'язком еліптичної крайової задачі \eqref{13f1}, \eqref{13f2}, праві частини якої задовольняють умови
\begin{gather*}
f\in H^{\eta_{1}\rho^{-2q}}_{\mathrm{loc}}(\Omega,\varnothing)\cap H^{\eta_{2}\rho^{-2q}}_{\mathrm{loc}}(U_{\delta},\Gamma),\\
g_{j}\in H^{\eta_{2}\rho^{-m_j-1/2}}(\Gamma),
\quad\mbox{при кожному}\quad j=1,...,q+\lambda
\end{gather*}
для деяких параметрів $\eta_{1},\eta_{2}\in\mathrm{RO}$ таких, що
$\sigma_{0}(\eta_{1})>m+1/2$, $\sigma_{0}(\eta_{2})>m+1/2$ і
\begin{equation*}
\int\limits_1^{\infty}t^{4q+n-1}\eta_{1}^{-2}(t)dt<\infty,\quad
\int\limits_1^{\infty}t^{2m+n-1}\eta_{2}^{-2}(t)dt<\infty.
\end{equation*}
Тоді цей розв'язок класичний.
\rm

\textbf{5. Випадок однорідного еліптичного рівняння.} Якщо в еліптичному рівнянні \eqref{13f1} права частина $f=0$ в області $\Omega$, то версії теорем 1--3 є правильними без обмеження $\sigma_{0}(\eta)>m+1/2$ на параметр $\eta\in\mathrm{RO}$. Наведемо, для прикладу, відповідну версію теореми 1.

Позначимо через $C^{\infty}(\overline{\Omega},A)$ лінійний простір усіх функцій $u\in C^{\infty}(\overline{\Omega})$, які задовольняють умову $Au=0$ в області $\Omega$. Окрім того, позначимо через $H^{\eta}(\Omega,A)$, де $\eta\in\mathrm{RO}$, лінійний простір усіх розподілів $u\in H^{\eta}(\Omega)$, які задовольняють цю умову (тепер вона розуміється у сенсі розподілів). Цей простір наділяємо нормою з $H^{\eta}(\Omega)$, відносно якої він є гільбертовим.

\textbf{Теорема 7.} \it Розглянемо відображення $(u,v_{1},...,v_{\lambda})\to(g_{1},...,g_{q+\lambda})$, де $u\in C^{\infty}(\overline{\Omega},A)$ та $v_{1},\ldots,v_{\lambda}\in C^{\infty}(\Gamma)$, а функції $g_{1}$,...,$g_{q+\lambda}$ означені за формулою \eqref{13f2}. Воно продовжується єдиним чином (за неперервністю) до обмеженого лінійного оператора
\begin{equation*}
\Lambda':H^{\eta}(\Omega,A)\oplus
\bigoplus_{k=1}^{\lambda}H^{\eta\rho^{r_{k}-1/2}}(\Gamma)\to
\bigoplus_{j=1}^{q+\lambda}H^{\eta\rho^{-m_{j}-1/2}}(\Gamma)
\end{equation*}
для довільного параметра $\eta\in\mathrm{RO}$. Цей оператор нетерів. Його скінченновимірне ядро лежить у просторі $C^{\infty}(\overline{\Omega},A)\times(C^{\infty}(\Gamma))^{\lambda}$ і разом зі скінченним індексом не залежить від $\eta$. \rm

За припущення, що функція $\eta$ правильно змінна на нескінченності, ця теорема встановлена недавно в [\ref{Anop19Dop2}, теорема 1].

\bigskip\bigskip

\noindent ЦИТОВАНА ЛІТЕРАТУРА

\begin{enumerate}

\item\label{LionsMagenes71}
Лионс Ж.-Л., Мадженес Э. Неоднородные граничные задачи и их приложения. Москва: Мир, 1971. 372 с.

\item\label{Lawruk63a}
Лаврук Б. О параметрических граничных задачах для эллиптических систем линейных дифференциальных уравнений. I. Построение сопряженных задач. \textit{Bull. Acad. Polon. Sci. S\'{e}r. Sci. Math. Astronom. Phys.} 1963. \textbf{11}, №~5. С. 257--267.

\item\label{KozlovMazyaRossmann97}
Kozlov V.A., Maz'ya V.G., Rossmann J. Elliptic boundary value problems in domains with point singularities. Providence: American Math. Soc.,  1997. ix+414 p.

\item\label{Roitberg99}
Roitberg Ya.A.  Elliptic boundary value problems in the spaces of distributions. Dordrecht: Kluwer Acad. Publisher, 1999. x+276~p.

\item\label{MikhailetsMurach13UMJ3}
Михайлец В.А., Мурач А.А. Расширенная соболевская шкала и эллиптические операторы. \textit{Укр. мат. журн.} 2013. \textbf{65}, №~3. С. 368--380.

\item\label{MikhailetsMurach14}
Mikhailets V.A., Murach A.A. H\"ormander spaces, interpolation, and elliptic problems. Berlin, Boston: De Gruyter, 2014. xii+297 p.

\item\label{Seneta85}
Сенета Е. Правильно меняющиеся функции. Москва: Наука, 1985. 144 с.

\item\label{BinghamGoldieTeugels89}
Bingham N.H., Goldie C.M., Teugels J.L. Regular variation. Cambridge: Cambridge Univ. Press, 1989. xix+494 p.

\item\label{Hermander65}
Хермандер Л. Линейные дифференциальные операторы с частными производными. Москва: Мир, 1965. 380 с.

\item\label{VolevichPaneah65}
Волевич Л.Р., Панеях Б.П. Некоторые пространства обобщенных функций и теоремы вложения. \textit{Успехи мат. наук.} 1965. \textbf{20}, № 1. С. 3--74.

\item\label{MikhailetsMurach15ResMath1}
Mikhailets V.A., Murach A.A. Interpolation Hilbert spaces between Sobolev spaces. \textit{Results Math.} 2015. \textbf{67}, №~1. P. 135--152.

\item\label{ChepuruhinaKasirenko17}
Касіренко Т.М., Чепурухіна І.С. Еліптичні за Лавруком задачі з крайовими операторами вищих порядків в уточненій соболєвській шкалі. \textit{Збірник праць Інституту математики НАН України}. 2017. \textbf{14}, № 3. С. 161--204.

\item\label{KasirenkoMurachUMJ3}
Касіренко Т.М., Мурач О.О. Еліптичні задачі з крайовими умовами високих порядків у просторах Хермандера. \textit{Укр. мат. журн.} 2017. \textbf{69}, № 11. С. 1486--1504.

\item\label{Chepuruhina15Coll2}
Чепурухіна І.С. Еліптичні крайові задачі за Б.~Лавруком у розширеній соболєвській шкалі. \textit{Збірник праць Інституту математики НАН України}. 2015. \textbf{12}, №~2. С. 338--374.

\item\label{Anop19Dop2}
Аноп А.В. Еліптичні за Лавруком крайові задачі для однорідних диференціальних рівнянь. \textit{Допов. Нац. акад. наук. Укр.} 2019. № 2. С. 3--11.

\end{enumerate}

\bigskip

\noindent REFERENCES

\begin{enumerate}

\item
Lions, J.-L. \& Magenes, E. (1972). Non-homogeneous boundary-value problems and applications, vol.~I. Berlin: Springer.

\item
Lawruk, B. (1963). Parametric boundary-value problems for elliptic systems of linear differential equations. I. Construction of conjugate problems. Bull. Acad. Polon. Sci. S\'{e}r. Sci. Math. Astronom. Phys., 11, No. 5, pp. 257–267 (in Russian).

\item
Kozlov, V.A., Maz'ya, V.G. \& Rossmann, J. (1997). Elliptic boundary value problems in domains with point singularities. Providence: Amer. Math. Soc.

\item
Roitberg, Ya.A. (1999). Elliptic boundary value problems in the spaces of distributions. Dordrecht: Kluwer Acad. Publ.

\item
Mikhailets, V.A. \& Murach, A.A. (2013). Extended Sobolev scale and elliptic operators. Ukrainian Math. J., 65, No. 3, pp. 435-447.

\item
Mikhailets, V.A. \& Murach, A.A. (2014). H\"ormander spaces, interpolation, and elliptic problems. Berlin, Boston: De Gruyter.

\item
Seneta, E. (1976). Regularly varying functions. Berlin: Springer.

\item
Bingham, N.H., Goldie, C.M. \& Teugels, J.L. (1989). Regular variation. Cambridge: Cambridge University Press.

\item
H\"ormander, L. (1963). Linear partial differential operators. Berlin: Springer.

\item
Volevich, L.R. \& Paneah, B.P. (1965). Certain spaces of generalized functions and embedding  theorems. Russian Math. Surveys, 20, No. 1, pp. 1-73.

\item
Mikhailets, V.A. \& Murach, A.A. (2015). Interpolation Hilbert spaces between Sobolev spaces. Results Math, 67, No. 1, pp. 135-152.

\item
Kasirenko, T.M. \& Chepurukhina, I.S. (2017). Elliptic problems in the sense of Lawruk with boundary operators of higher orders in refined Sobolev scale. Zbirnyk Prats Institutu Matematyky NAN Ukrainy, 14, No. 3, pp. 161-204 (in Ukrainian).

\item
Kasirenko, T.M. \& Murach, O.O. (2018). Elliptic problems with boundary conditions of higher orders in H\"ormander spaces. Ukrainian Math. J., 69, No. 11, pp. 1727-1748.

\item
Chepurukhina, I.S. (2015). Elliptic boundary-value problems in the sense of Lawruk in an extended Sobolev scale. Zbirnyk Prats Institutu Matematyky NAN Ukrainy, 12, No. 2, pp. 338-374 (in Ukrainian).

\item
Anop, A.V. (2019). Lawruk elliptic boundary-value problems for homogeneous differential equations. Dopov. Nac. akad. nauk Ukr., No 2, pp. 3-11 (in Ukrainian).

\end{enumerate}

\end{document}